 \documentclass[final,review,3p,times]{elsarticle}

\usepackage{amssymb}

\usepackage{graphicx}

\usepackage{url}

\journal{}

\begin{document}

\begin{frontmatter}


\newtheorem{theorem}{Theorem}[section]
\newtheorem{lemma}{Lemma}[section]
\newtheorem{algorithm}{Algorithm}[section]
\newtheorem{remark}{Remark}[section]
\newtheorem{assumption}{Assumption}[section]

\newtheorem{question}[theorem]{Question}
\newcommand{\proof}{\nr \bf Proof: \hspace{0mm} \rm }

\title{An Improved Sequential Quadratic Programming Algorithm
 for Solving General Nonlinear Programming Problems
 \tnoteref{s1}}
 \tnotetext[s1]{Guo and Bai's research is supported by the National
Natural Science Foundation of China (Grant No. 11071158), Jian's
research is supported by the National Natural Science Foundation of
China ( Grant No. 11271086) and the Natural Science Foundation of
Guangxi (Grant No. 2011GXNSFD018002) as well as Innovation Group of
Talents Highland of Guangxi higher School.}


\author[shu]{Chuan-Hao Guo\corref{cor1}}
\ead{guo-ch@live.cn}

\author[shu]{Yan-Qin Bai\corref{cor2}}
\ead{yqbai@shu.edu.cn}

\author[gxu]{Jin-Bao Jian\corref{cor3}}
\ead{jianjb@gxu.edu.cn}

\cortext[cor2]{Corresponding author}

\address[shu]{Department of Mathematics, Shanghai University,
Shanghai 200444, China}

\address[gxu]{Department of Mathematics and Information Science, Guangxi University, Nanning 530004, Guangxi, China;
Headmaster's Office, Yulin Normal University, Yulin 537000, Guangxi, China}

\begin{abstract}
In this paper, a class of general nonlinear programming problems
with inequality and equality constraints is discussed. Firstly, the
original problem is transformed into an associated simpler
equivalent problem with only inequality constraints. Then, inspired
by the ideals of sequential quadratic programming (SQP) method and
the method of system of linear equations (SLE), a new type of SQP
algorithm for solving the original problem is proposed. At each
iteration, the search direction is generated by the combination of
two directions, which are obtained by solving an always feasible
quadratic programming (QP) subproblem and a SLE, respectively.
Moreover, in order to overcome the Maratos effect, the higher-order
correction direction is obtained by solving another SLE. The two
SLEs have the same coefficient matrices, and we only need to solve
the one of them after a finite number of iterations. By a new line
search technique, the proposed algorithm possesses global and
superlinear convergence under some suitable assumptions without the
strict complementarity. Finally, some comparative numerical results
are reported to show that the proposed algorithm is effective and
promising.
\end{abstract}

\end{frontmatter}

\section{Introduction}\label{sec1}
In this paper we consider the following nonlinear programming
problem
\begin{equation}\label{problem1}
\begin{array}{lll}
&\min\ \ f_0(x)\\
&{\rm s.t.} \ \ \ f_i(x)\leq 0,\ \ i\in I_1:=\{1,2,\dots,m'\},\\
&\ \ \ \ \ \ \  f_i(x)=0,\ \ i\in I_2:=\{m'+1,m'+2,\dots,m\},
 \end{array}
 \end{equation}
  where $f_{i}\ : R^{n} \rightarrow R \ (i\in\{0\}\cup I_1\cup I_2)$
  are smooth functions. The feasible set and gradients
  of problem (\ref{problem1}) are denoted as follows:
$$\Omega:=\{x\in R^n:\ f_{i}(x)\leq 0, i\in I_1;\ f_{i}(x)=0,
i\in I_2\},\ {\rm{and}}\ g_i(x):=\nabla f_i(x),\ i\in\{0\}\cup
I_1\cup I_2.$$

Sequential quadratic programming (SQP) algorithms have been widely
studied by many authors during the past several decades, e.g., Refs.
\cite{s1998i,qy2001,jthz2008,zzg2010,gz2011,cllpy2012,zll2012}, and
have been proved highly effective for solving problem
(\ref{problem1}). SQP algorithms generate iteratively the main
search directions by solving the standard quadratic programming (QP)
subproblem
$$\begin{array}{lll}
&\min\ \ g_0(x)^Td+\frac{1}{2}d^THd\\
&{\rm s.t.} \ \ \ f_i(x)+g_i(x)^Td\leq 0,\ \ i\in I_1,\\
&\ \ \ \ \ \ \ f_i(x)+g_i(x)^Td=0,\ \ i\in I_2,
 \end{array}
$$ where $H\in R^{n\times n}$ is a symmetric positive definite matrix.
Then one performs a line search which is a one dimensional
minimization problem to determine a steplength, and obtain the next
iteration point.

SQP algorithms may fail since the equality constraints of QP
subproblem are hard to be satisfied in the process of iteration.
Mayne and Polak \cite{mp1976} propose a new way for overcoming this
difficulty. In their scheme, they consider the following related
family of simpler problem
\begin{equation}\label{problem2}\ \ \ \begin{array}{lll}
&\min\ \ F_c(x):= f_0(x)-c\sum\limits_{i\in I_2}f_i(x)\\
&{\rm s.t.} \ \ \ f_i(x)\leq 0,\ \ i\in I,
 \end{array}
 \end{equation}
 where $I:=I_1\cup I_2$ and parameter $c>0$. Especially, $F_c(x)=f_0(x)$ if $I_2=\emptyset$. We denote
 the feasible set of problem (\ref{problem2}) by
 $$\Omega^+:=\{x\in R^n:\ f_i(x)\leq 0,\ i\in I\}.$$
 Moreover, they prove that the original problem (\ref{problem1}) is equivalent to problem (\ref{problem2}) when $c$ is
 sufficiently large but finite. Note that problem (\ref{problem2}) only has inequality
 constraints, so the corresponding QP subproblem has not equality
 constraints, and SQP algorithms will be always successful under some
 suitable conditions. More advantages and
 further applications of this technique can be seen in
 \cite{lt1996,jgy2009,h1986}.

Recently, Guo propose an algorithm for solving problem
(\ref{problem1}) with $I_2=\emptyset$ in \cite{g2011}. In this
algorithm, the initial iteration point can be chosen arbitrarily.
The main search direction is obtained by solving one QP subproblem
and one (or two) system(s) of linear equations. The algorithm
possesses global and superlinear convergence under some suitable
assumptions without the strict complementarity. Furthermore, some
comparative numerical results are reported to show that the
algorithm is effective.

Inspired by the ideas in \cite{mp1976,g2011}, we propose a new SQP
algorithm for solving  problem (\ref{problem1}). First,  problem
(\ref{problem1}) is equivalently transformed into  problem
(\ref{problem2}) (see Lemma \ref{lemma2.3}). In order to overcome
the inconstant of QP subproblem,  we consider a modified QP
subproblem
\begin{equation}\label{problem3}\begin{array}{lll}
&\min\ \ \nabla F_c(x)^Td+\frac{1}{2}d^THd\\
&{\rm s.t.} \ \ \ f_i(x)+g_i(x)^Td\leq\varphi(x), &i\in I^+(x),\\
& \ \ \ \ \ \ \ f_i(x)+g_i(x)^Td\leq 0, &i\in I^-(x),
 \end{array}
\end{equation} where $\varphi(x):=\max\{0,f_i(x), i\in I\},\ I^+(x):=\{i\in I:\
f_i(x)>0\},\ I^-(x):=\{i\in I:\ f_i(x)\leq 0\}$. QP subproblem
(\ref{problem3}) has the following advantages:

$\bullet$ subproblem (\ref{problem3}) always has a feasible solution
$d = 0$.

$\bullet$ subproblem (\ref{problem3}) is a strictly convex program
if
 $H$ is positive definite, so it always has a unique solution.

$\bullet$ $d$ is a solution of subproblem (\ref{problem3}) if and
only if it is a KKT point of subproblem (\ref{problem3}).

In order to get the global convergence of the algorithm, the search
direction is generated by the combination of two directions, which
are obtained by solving QP subproblem (\ref{problem3}) and a system
of linear equations, respectively. For overcoming the Maratos effect
\cite{m1978}, the higher-order correction direction is generated by
solving another system of linear equations. The two systems of
linear equations have the same coefficient matrices. The superlinear
convergence is derived under the strong second-order sufficient
conditions (SSOSC) without the strict complementarity. Moreover, for
further comparing the performance of the method of strongly
sub-feasible directions (MSSFD) \cite{jthz2008,jgy2009,g2011} with
the method of quasi-strongly sub-feasible directions (MQSSFD)
\cite{jcg2011}, the technical of MQSSFD is adopted in our new
algorithm.  Finally, some comparative numerical results are reported
to show that our new algorithm is promising. The main features of
the proposed algorithm are summarized as follows:

$\bullet$ the initial iteration point is arbitrary, and the number
of constraints satisfying constraint condition is monotone
nondecreasing.

$\bullet$ the objective function of problem (\ref{problem2}) is used
directly as the merit function.

$\bullet$ the parameter $c$ is adjusted automatically only for a
finite number of times (see Lemma \ref{lemma3.2}).

$\bullet$ at each iteration, the search direction is generated by a
combination of two directions, which are obtained by solving an
always feasible QP subproblem and a system of linear equations,
respectively.

$\bullet$ after finite iterations, the iteration points always lie
in $\Omega^+$.

$\bullet$ under SSOSC without the strict complementarity, the
proposed algorithm possesses global and superlinear convergence.

The paper is organized into six sections. In Section \ref{sec2}, our
new algorithm and its properties are presented. In Sections
\ref{sec3} and \ref{sec4}, we show that the proposed algorithm
possesses global and superlinear convergence, respectively. In
Section \ref{sec5}, some comparative numerical results are reported
to show that the proposed algorithm is effective and promising. Some
conclusions about the proposed algorithm are given in Section
\ref{sec6}.

Throughout the paper we use the following notations for a point
$x\in R^n$ and an index subset $J\subseteq I$
$$\left\{\begin{array}{ll}f_J(x):=\{f_i(x),\ i\in J\},\ g_J(x):=\{g_i(x),\ i\in
J\},\ \bar{f}_i(x):=f_i(x),\ i\in I^-(x),\
\bar{f}_i(x):=f_i(x)-\varphi(x),\ i\in I^+(x),\\
I_1(x):=\{i\in I_1:\ \bar{f}_i(x)=0\},\ I(x):=I_1(x)\cup I_2,\
I_0(x):=\{i\in I:\ \bar{f}_{i}(x)=0\}.
 \end{array}\right.$$

\section{Description of algorithm}\label{sec2}
In this section, we start by giving some basic assumptions for
problem (\ref{problem1}).
\begin{assumption}\label{assumption 2.1}{\it (i) The functions $f_i(x)\ (i\in\{0\}\cup I)$ are all
continuously differentiable.\\
(ii) The gradient vectors $\{g_{i}(x):\ i\in {I}(x)\}$ are linearly
independent for each $x\in R^n$.}
\end{assumption}

To update the parameter $c$ in  problem (\ref{problem2}), the matrices
$N(x),\ D(x)$ and multiplier vector $\pi(x)$ are defined as follows:
\begin{equation}\label{2.1}\begin{array}{ll}N(x)=(g_i(x),i\in I),\
\pi(x)=-(N(x)^TN(x)+D(x))^{-1}N(x)^Tg_0(x),\\
D(x)={\rm{diag}}(D_i(x), i\in
I),\ D_i(x)=\left\{\begin{array}{ll}|\bar{f}_i(x)|^p,\ &i\in I_1,\\
0,\ \ \ \ \ &i\in I_2,\end{array}\right.\\
\end{array}\end{equation}
where $p$ is a positive parameter.

Note that $D_i(x)>0$ for all $i\in I\backslash I(x)$ in (\ref{2.1}).
By Assumption \ref{assumption 2.1}(ii), the following lemma holds
immediately.
\begin{lemma}\label{lemma2.2} Suppose that Assumption \ref{assumption 2.1} holds.
Then $(N(x)^TN(x)+D(x))$ is nonsingular and positive definite for
all $x\in R^n$.
\end{lemma}

By Lemma \ref{lemma2.2} and (\ref{2.1}), the relationship between
problems (\ref{problem1}) and (\ref{problem2}) is shown in the
following Lemma \ref{lemma2.3}. Its proof can be referred to the one
of Lemma 2.1(v) in \cite{jgy2009}.
\begin{lemma}\label{lemma2.3} If $c>|\pi_i(x)|$ for all $i\in
I_2$, then $(x, \mu)$ is the KKT point of  problem (\ref{problem1})
if and only if $(x, \lambda)$ is the KKT point of problem
(\ref{problem2}), where $\mu$ and $\lambda$ satisfy
\begin{equation}\label{2.2}\mu_i=\lambda_i,\ i\in I_1,\ \ \ \
\mu_i=\lambda_i-c,\ i\in I_2.\end{equation}
\end{lemma}

For the iteration point $x^k$ and the parameter $c_k$ of problem
(\ref{problem2}), QP subproblem (\ref{problem3}) can be simplified
as follows by the above notations
\begin{equation}\label{problem2.8}\ \ \ \begin{array}{lll}
&\min\ \ \nabla F_{c_k}(x^k)^Td+\frac{1}{2}d^TH_kd\\
&{\rm s.t.} \ \ \ \bar{f}_i(x^k)+g_i(x^k)^Td\leq 0,\ \ i\in I.
 \end{array}
 \end{equation}
When $H_k$ is positive definite, $d_0^k$ is a solution of subproblem
(\ref{problem2.8}) if and only if there exists a corresponding KKT
multiplier vector $\lambda^k$ such that
\begin{equation}\label{2.9}\left\{\begin{array}{ll}\nabla
F_{c_k}(x^k)+H_kd_0^k+\sum\limits_{i\in
I}\lambda_i^kg_i(x^k)=0,\\
\bar{f}_i(x^k)+g_i(x^k)^Td_0^k\leq 0,\ \lambda_i^k\geq 0,\
\lambda_i^k(\bar{f}_i(x^k+g_i(x^k)^Td_0^k))=0,\ \forall\ i\in
I.\end{array}\right.\end{equation} Since $d=0$ is a feasible
solution of subproblem (\ref{problem2.8}) and $H_k$ is positive
definite, it follows that
$$\nabla F_{c_k}(x^k)^Td_0^k+\frac{1}{2}(d_0^k)^TH_kd_0^k\leq0\Rightarrow \nabla F_{c_k}(x^k)^Td_0^k\leq 0,$$
i.e., $d_0^k$ is a descent direction of $F_{c_k}(x^k)$ at the
iteration point $x^k$.

Due to (\ref{2.9}) and Lemma \ref{lemma2.3}, the following lemma
holds immediately.
\begin{lemma}\label{lemma2.4}
If $(d_0^k, \varphi(x^k))=(0, 0)$, then $x^k$ is the KKT point of
problem (\ref{problem2}). Furthermore, if $c_k>|\mu_i^k|$ for all
$i\in I_2$, then $x^k$ is the KKT point of  problem
(\ref{problem1}).
\end{lemma}

Again from (\ref{2.9}), it follows that $d_0^k$ may not be a
feasible direction of problem (\ref{problem2}) at the feasible
iteration point $x^k\in\Omega^+$. So a suitable strategy must be
carried out to generate a feasible direction. Here, taking into
account that $x^k$ may be infeasible, we introduce the system of
linear equations to get a unique solution $(d_1^k, h_1^k)$
 \begin{equation}\label{2.11} \Gamma_k\left(\begin{array}{c}d\\h
\end{array}\right)\triangleq\left(\begin{array}{c}H_k\ \ \ N_k\\N_k^T\ \
-Q^k\end{array}\right)\left(\begin{array}{c}d\\h
\end{array}\right)=\left(\begin{array}{c}0\\-(||d_0^k||+\varphi(x^k)^\sigma)\varpi\end{array}\right),
\end{equation}
 where $0\in
R^n$, $\varpi=(1,1,\ldots,1)^T\in R^m$, $\sigma\in(0,1)$ and
 \begin{equation}\label{2.12}\begin{array}{ll}N_k=N(x^k)=(g_i(x^k),i\in I),\ \
Q^k={\rm
diag}(Q_{i}^k=|\bar{f}_i(x^k)|(|\bar{f}_i(x^k)+g_i(x^k)^Td_0^k|+||d_0^k||),\
i\in I).\end{array}\end{equation}
 Then we consider the convex combination of $d_0^k$ and $d_1^k$
\begin{equation}\label{2.13}\hat{d}^k=(1-\beta_k)d_0^k+\beta_kd_1^k,\end{equation}
 where
$\beta_k$ is the maximal value of $\beta\in[0, 1]$ satisfying
 \begin{equation}\label{2.14}\nabla
F_{c_k}(x^k)^T\hat{d}^k\leq\theta \nabla
F_{c_k}(x^k)^Td_0^k+\varphi(x^k)^\theta.\end{equation} Moreover,
(\ref{2.14}) further implies that $\beta_k$ is the optimal solution
of linear programming
$$\begin{array}{ll} \max\ \ \ \ \ \beta \\ \ {\rm s.t.} \ \ \
\beta \nabla F_{c_k}(x^k)^Td_1^k+(1-\beta)\nabla
F_{c_k}(x^k)^Td_0^k\leq \theta\nabla
 F_{c_k}(x^k)^Td_0^k+\varphi(x^k)^\theta,\\
\ \ \ \ \ \ \ \ 0\leq\beta\leq1,\end{array} $$
 where the positive
parameter $\theta<\sigma$. It is obviously that (\ref{2.14}) holds
for $\beta_k=0$, since $\varphi(x^k)\geq 0$. The above linear
programming further implies that $\beta_k>0$, otherwise, $x^k$ is
the KKT point of problem (\ref{problem1}) (see Lemma
\ref{lemma2.4}).

The next lemma shows the solvability of (\ref{2.11}). Its proof is
elementary in view of $\{i\in I:\ Q_i^k=0\}\subseteq
I_0(x^k)\subseteq I(x^k)$.
\begin{lemma}\label{lemma2.5} Suppose that Assumption \ref{assumption
2.1} holds and $H_k$ is positive definite. Then $\Gamma_k$ defined
in (\ref{2.11}) is nonsingular and (\ref{2.11}) has a unique
solution.
\end{lemma}

\begin{lemma}\label{lemma2.6} Suppose that
Assumption \ref{assumption 2.1} holds. Then\\
(i)
$\nabla F_{c_k}(x^k)^T\hat{d}^k\leq-\frac{1}{2}\theta(d_0^k)^TH_kd_0^k+\varphi(x^k)^\theta$.\\
(ii) $g_i(x^k)^T\hat{d}^k\leq
-\beta_k(||d_0^k||+\varphi(x^k)^\sigma),\ \forall\ i\in{I_0}(x^k).$
\end{lemma}
{\bf proof:} (i) Since $d=0$ is a feasible solution of  subproblem
(\ref{problem2.8}), and from (\ref{2.14}), it holds that
 $$\begin{array}{ll}\nabla F_{c_k}(x^k)^T\hat{d}^k\leq\theta\nabla F_{c_k}(x^k)^T{d}_0^k+\varphi(x^k)^\theta
 \leq-\frac{1}{2}\theta(d_0^k)^TH_kd_0^k+\varphi(x^k)^\theta.\end{array}$$

(ii) First, from (\ref{2.9}) and (\ref{2.11}), it follows that
$$g_i(x^k)^Td_0^k\leq 0,\
g_i(x^k)^T{d}_1^k=-||d_0^k||-\varphi(x^k)^\sigma,\ \forall\ i\in
I_0(x^k).$$ Then we obtain that by (\ref{2.13})
$$g_i(x^k)^T\hat{d}^k \leq-\beta_k(||d_0^k||+\varphi(x^k)^\sigma)\
{\rm for\ all} \ i\in I_0(x^k).$$\qed

By Lemma \ref{lemma2.6}, we know that $\hat{d}^k$ is an improved
direction. In order to overcome the Maratos effect and avoid the
strict complementarity condition as well as reduce the computational
cost, a suitable higher-order correction direction should be
introduced by an appropriate approach. Here, we introduce the
following system of linear equations to yield the higher-order
correction direction $d_2^k$
 \begin{equation}\label{2.15}
\Gamma_k\left(\begin{array}{c}d\\h
\end{array}\right)=\left(\begin{array}{c}0\\-(||d_0^k||^\tau+\varphi(x^k)^\sigma)\varpi-\digamma(x^k+d_0^k)\end{array}\right),
\end{equation}
 where
$\tau\in(2,3)$ and
 \begin{equation}\label{2.16}\digamma(x^k+d_0^k)=(f_i(x^k+d_0^k)-f_i(x^k)-g_i(x^k)^Td_0^k,\ i\in
I).\end{equation} Note that the term $\varphi(x^k)$ is introduced in
our paper, the relationship between $d_2^k$ and $d_0^k$ will be
different from the traditional form  $||d_2^k||=O(||d_0^k||^2)$
\cite{hc1986,pt1987,pt1993}, the details can be seen in Lemma
\ref{lemma4.4}.

We are now ready to present our algorithm for solving problem
(\ref{problem1}) as follows.
\begin{algorithm}\label{algorithm 2.7}\end{algorithm}

Parameters: $p,\epsilon,\gamma,\gamma_0>0,\ c_{-1}>0,\ \rho>1,\
0<\theta<\sigma,\ \sigma,\eta,\alpha,\hat{\alpha}\in(0,1),\
\tau\in(2,3).$

Data: $x^0\in R^n$, a symmetric positive definite matrix $H_0\in
R^{n\times n}$, and $k:=0$.

\textbf{Step 1.} Update parameter $c_k$: Compute $c_k\
(k=1,2,\dots)$ by
\begin{equation}\label{2.17} c_k=\left\{\begin{array}{ll}\max\{s_k,\ c_{k-1}+\gamma\}, & {\rm if}\  s_k>c_{k-1},\\
  c_{k-1}, & {\rm if}\  s_k\leq c_{k-1},\\
  \end{array}\right.s_k=\max\{|\pi_{i}(x^k)|,\ i\in I_{2}\}+\gamma_0.
 \end{equation}

\textbf{Step 2.} Solve QP subproblem: Solve QP subproblem
(\ref{problem2.8}) to get a solution. If $(d_0^k,
\varphi(x^k))=(0,0)$, then $x^k$ is the KKT point of problem
(\ref{problem1}) and stop; otherwise, go to Step 3.

\textbf{Step 3.} Solve system of linear equations: Solve
(\ref{2.15}) to get a solution $(d_2^k,h_2^k)$, and let
$d^k=d_0^k+d_2^k$.

\textbf{Step 4.} Let $t=1$. (a) If
 \begin{equation}\label{2.18}\left\{\begin{array}{ll}
 F_{c_k}(x^k+td^k)\leq
F_{c_k}(x^k)+\alpha t \nabla F_{c_k}(x^k)^Td_0^k+\rho(1-\alpha)t\varphi(x^k)^\theta,\\
  f_i(x^k+td^k)\leq\max\{0,\varphi(x^k)-\alpha t(||d_0^k||^\tau+\varphi(x^k)^\sigma)\},\ \ i\in
 I,\\
  |I^-(x^k+t d^k)|\geq |I^-(x^k)|,\\
 \end{array}\right.
 \end{equation}
 is satisfied, then let $t_k=t$, and go to Step 7; otherwise, go to (b).

 $\ \ \ \ \ \ \ \ \ \ \ \ \ \ \ \ \ {\rm{(b)}}\ {\rm{Let}}\
 t:=\frac{1}{2}t.$\\
If $t<\epsilon$, then go to Step 5; otherwise, repeat (a).

\textbf{Step 5.} Solve system of linear equations: Solve
(\ref{2.11}) to get a  solution $(d_1^k, h_1^k)$, and compute
$\hat{d}^k$ by (\ref{2.13}) and (\ref{2.14}).

\textbf{Step 6.} Compute steplength $t_k$ be the first member of the
sequence $\{1,\eta,\eta^2,\ldots\}$ such that
 \begin{equation}\label{2.19}
\left\{\begin{array}{ll} F_{c_k}(x^k+t\hat{d}^k)\leq
F_{c_k}(x^k)+\hat{\alpha} t
\nabla F_{c_k}(x^k)^T\hat{d}^k+\rho(1-\hat{\alpha})t\varphi(x^k)^\theta,\\
f_i(x^k+t\hat{d}^k)\leq\max\{0,\ \varphi(x^k)-\hat{\alpha}
t{\beta}_k(||d_0^k||+\varphi(x^k)^\sigma)\},\ \ i\in
 I,\\
 |I^-(x^k+t \hat{d}^k)|\geq |I^-(x^k)|,\end{array}\right.\end{equation}
 and let $d^k=\hat{d}^k$.

\textbf{Step 7.} Compute a new symmetric positive definite matrix
$H_{k+1}$ by some techniques, set $x^{k+1}:=x^k+t_k d^k,\ k:=k+1$,
and go to Step 1.
\begin{lemma}\label{lemma2.9} Suppose that Assumption \ref{assumption 2.1} holds.
If Algorithm \ref{algorithm 2.7} does not stop at Step 2, i.e.,
$(d_0^k, \varphi(x^k))\neq(0,0)$, then the line search (\ref{2.19})
can be terminated after a finite number of iterations.
\end{lemma}
In fact, the above Lemma \ref{lemma2.9} shows that Algorithm
\ref{algorithm 2.7} is well defined. Moreover, $\hat{d}^k$ is an
improved direction in a sense by Lemma \ref{lemma2.6}. Therefore, if
Algorithm \ref{algorithm 2.7} does not stop at Step 2, i.e.,
$(d_0^k,\varphi(x^k))\neq (0,0)$, we can always get the next
iteration point $x^{k+1}$ from the current iteration point $x^k$
according to Lemma \ref{lemma2.9}. Furthermore, from the mechanism
of Algorithm \ref{algorithm 2.7}, the following lemma holds
obviously.
\begin{lemma}\label{lemma2.10} Suppose that Assumption
\ref{assumption 2.1} holds.\\
(i) If there exists an index $k_0$ such that $x^{k_0}\in \Omega^+$,
then $x^k\in \Omega^+$ and
$F_{c_k}(x^{k+1})\leq F_{c_k}(x^k)$ for all $k\geq k_0$.\\
(ii) If $x^k\not\in \Omega^+$ and $x^{k+1}\not\in \Omega^+$, then
$\varphi(x^{k+1})<\varphi(x^k)$.\\
(iii) The subsets $I^-(x^k)$ and $I^+(x^k)$ can be fixed, i.e.,
$I^-(x^k)\equiv I^-$ and $I^+(x^k)\equiv I^+$ for $k$ large enough.
\end{lemma}
\section{Global convergence}\label{sec3}
In this section, we establish the global convergence of Algorithm
\ref{algorithm 2.7}. If Algorithm \ref{algorithm 2.7} stops at
$x^k$, it follows that $x^k$ is the KKT point of problem
(\ref{problem1}).  Now, we assume that Algorithm \ref{algorithm 2.7}
produces an infinite sequence $\{x^k\}$ of iteration points, and
prove that each accumulation point $x^*$ of $\{x^k\}$ is the KKT
point of problem (\ref{problem1}) under some suitable assumptions.
For this purpose, the following assumption is necessary.
\begin{assumption}\label{assumption 3.1} (i) The sequence $\{x^k\}$ is
bounded.\\
(ii) There exist positive constants $a$ and $b$ such that
\begin{equation}\label{3.1}a||d||^2\leq d^TH_kd\leq b||d||^2,\ \ \forall\ d\in
R^n,\ \forall\ k.\end{equation}
\end{assumption}

Denote the active set for QP subproblem (\ref{problem2.8}) by
$$L(x^k)\equiv\{i\in I:\ \bar{f}_i(x^k)+g_i(x^k)^Td_0^k=0\}.$$ Suppose
that $x^*$ is a given accumulation point of $\{x^k\}$. In view of
$I^+(x^k),\ I^-(x^k)$ and $L(x^k)$ are subsets of the finite set
$I$, by Lemma \ref{lemma2.10}(iii), we can assume that there exists
an infinite index set $K$ such that
 \begin{equation}\label{3.2}x^k\rightarrow x^*,\
I^-(x^k)\equiv I^-,\ I^+(x^k)\equiv I^+,\ L(x^k)\equiv L,\
\varphi(x^k)\rightarrow\varphi(x^*),\ \forall\ k\in K.\end{equation}

\begin{lemma}\label{lemma3.2}
Suppose that Assumptions \ref{assumption 2.1} and \ref{assumption
3.1}(i) hold. Then there exists an index $k_1>0$ such that
$c_k=c_{k_1}\triangleq c$ for all $k\geq k_1$.
\end{lemma}

The detailed proof of this lemma can be found in \cite{jgy2009}. Due
to Lemma \ref{lemma3.2}, we assume that $c_k\equiv c$ for all $k$ in
the rest of this paper. The results given in the following lemma are
very important in the subsequent analysis.
\begin{lemma}\label{lemma3.3} Suppose that Assumptions \ref{assumption 2.1} and \ref{assumption 3.1} hold. Then\\
(i) The sequence $\{d_0^k\}_{k=1}^\infty$ is bounded.\\
(ii) There exists a constant $r_0>0$ such that
$||\Gamma_k^{-1}||\leq r_0$
for all $k$.\\
(iii) The sequence $\{d_1^k\}_{k=1}^\infty$,
$\{d_2^k\}_{k=1}^\infty$, $\{\hat{d}^k\}_{k=1}^\infty$ and
$\{h_2^k\}_{k=1}^\infty$ are all bounded.
\end{lemma}
{\bf proof:}\ \ (i) Due to the fact that $d=0$ is a feasible
solution of subproblem (\ref{problem2.8}) and $d_0^k$ is an optimal
solution, we have \begin{equation}\label{3.3}\nabla
F_c(x^k)^Td_0^k+\frac{1}{2}(d_0^k)^TH_kd_0^k\leq 0.\end{equation} By
Assumption \ref{assumption 3.1} and the continuity of $\nabla
F_c(x^k)$, there exists a constant $\bar{c}>0$ such that $||\nabla
F_c(x^k)||\leq\bar{c}$ for all $k$. Combining (\ref{3.3}) with
(\ref{3.1}), we get
$$-\bar{c}||d_0^k||+\frac{1}{2}a||d_0^k||^2\leq 0,$$
which implies that $\{d_0^k\}$ is bounded for all $k$.

(ii) Suppose by contradiction that there exists an infinite index
set $K$ such that
 \begin{equation}\label{3.4}||\Gamma_k^{-1}||\rightarrow\infty,
\ k\in K.\end{equation} Without loss of generality, we assume that
there exists an infinite index set $K'\subseteq K$ such that
$$\begin{array}{ll}x^k\rightarrow x^*,\ \Gamma_k\rightarrow\Gamma_*=\left(\begin{array}{ll}H_*\ &N_*\\
N_*^T \ &-Q^*\end{array}\right),\ d_0^k\rightarrow d_0^*,\
\varphi(x^k)\rightarrow\varphi(x^*),\ I_0(x^k)\rightarrow I_0(x^*),
\ \forall\ k\in K',\end{array}$$
 where
$$\begin{array}{ll}N_*=(g_{i}(x^*),\ i\in I),\
Q^*={\rm
diag}(Q_{i}^*=|\bar{f}_i(x^*)|(|\bar{f}_i(x^*)+g_i(x^*)^Td_0^*|+||d_0^*||),\
i\in I).\end{array}$$
 It holds that $Q_i^*\geq 0$ for all $i\in I$, and $Q_i^*>0$ for all $i\in I\backslash I_0(x^*)$.
Similar to the proof of Lemma \ref{lemma2.5}, we can conclude that
$\Gamma_*$ is nonsingular. So
$||\Gamma_k^{-1}||\rightarrow||\Gamma_*^{-1}||$ for $k\in K'$, which
contradicts (\ref{3.4}), thus the conclusion (ii) holds.

(iii) Taking into account (\ref{2.11}), (\ref{2.15}) and
(\ref{2.13}), we can obtain the boundedness of
$\{d_1^k\}_{k=1}^\infty$, $\{d_2^k\}_{k=1}^\infty$,
$\{\hat{d}^k\}_{k=1}^\infty$ and $\{h_2^k\}_{k=1}^\infty$ by
employing the result of parts (i) and (ii).\qed

Similar to the analysis of Lemma 3.3 in \cite{g2011}, we can obtain
the following results.
\begin{lemma}\label{lemma3.4}
Suppose that Assumptions \ref{assumption 2.1} and \ref{assumption
3.1} are satisfied. Then\\
(i) $\lim\limits_{k\rightarrow\infty}(d_0^k, \varphi(x^k))=(0,0)$,
$\lim\limits_{k\rightarrow\infty}d_1^k=\lim\limits_{k\rightarrow\infty}d_2^k=
\lim\limits_{k\rightarrow\infty}\hat{d}^k=0$ and $\lim\limits_{k\rightarrow\infty}h_2^k=0$.\\
(ii) $\lim\limits_{k\rightarrow\infty}||x^{k+1}-x^k||=0$.
\end{lemma}
\begin{theorem}\label{theorem3.5} Suppose that Assumptions \ref{assumption 2.1} and \ref{assumption 3.1} hold.
Then Algorithm \ref{algorithm 2.7} either stops at the KKT point
$x^k$ of  problem {(\ref{problem1})} after a finite number of
iterations or generates an infinite sequence $\{x^k\}$ of points
such that each accumulation point $x^*$ of $\{x^k\}$ is the KKT
point of  problem {(\ref{problem1})}. Furthermore, there exists an
index set $K$ such that $\{(x^k,\lambda^k): k\in K\}$ and
$\{(x^k,\mu^k): k\in K\}$ converge to the KKT pair $(x^*,\lambda^*)$
of problem (\ref{problem3}) and the KKT pair $(x^*,\mu^*)$ of
problem (\ref{problem1}), respectively, where
$\lambda^k=(\lambda^k_L,0_{I\backslash L})$ and
$\mu^k=(\mu^k_L,0_{I\backslash L})$.
\end{theorem}
{\bf proof:}\ \ By Lemma \ref{lemma2.10}(iii), we assume without
loss of generality that there exists an infinite subset $K$ such
that (\ref{3.2}) holds. Let matrix $A_k=(g_i(x^k),\ i\in L)$. From
Lemma \ref{lemma3.4}(i), it follows that $L\subseteq I_0(x^*)=\{i\in
I:\ \bar{f}_i(x^*)=0\}$, which together with Assumption
\ref{assumption 2.1} shows that $A_k^TA_k$ is nonsingular for $k\in
K$ large enough, since $A_k\stackrel{K}{\longrightarrow}A_*$, where
$A_*\triangleq(g_j(x^*),\ i\in L).$

 From (\ref{2.9}) and Lemma \ref{lemma3.4}(i), we have for $k\in K$ large enough,
$$\lambda_{L}^k=-(A_k^TA_k)^{-1}A_k^T(\nabla F_{c}(x^k)+H_kd_0^k)\rightarrow-(A_*^TA_*)^{-1}A_*^T
\nabla F_{c}(x^*)\triangleq\lambda_{L}^*.$$
 Denote the
multiplier vector $\lambda^*=(\lambda_{L}^*,0_{I\setminus L})$, then
$\lim\limits_{k\in K}\lambda^k=\lambda^*.$ Passing to the limit
$k\in K$ ($k\rightarrow\infty$) in (\ref{2.9}), it follows that
$$\nabla F_{c}(x^*)+N_*\lambda^*=0,\ f_i(x^*)\leq 0,\ \lambda_i^*\geq 0,\ f_i(x^*)\lambda_i^*=0,\ i\in I,$$
which shows that $(x^*,\lambda^*)$ is the KKT pair of  problem
(\ref{problem3}). By the definition of $c_k$ and Lemma
\ref{lemma3.2}, we have $c>\max\{\pi_i(x^*): i\in I_2\}$. So from
Lemma \ref{lemma2.3}, we can conclude that $(x^*, \mu^*)$ is the KKT
pair of  problem (\ref{problem1}) with $\mu_i^*=\lambda_i^*,\ i\in
L\backslash I_2;\ \mu_i^*=\lambda_i^*-c,\ i\in I_2;\ \mu_i^*=0,\
i\in I\backslash L$. Obviously, $\lim\limits_{k\in
K}(x^k,\lambda^k)=(x^*,\lambda^*)$ and $\lim\limits_{k\in
K}(x^k,\mu^k)=(x^*,\mu^*)$. The proof is completed.\qed

\section{Rate of convergence}\label{sec4}
In this section we further discuss the strong and superlinear
convergence of Algorithm \ref{algorithm 2.7}. For these purposes, we
make the following assumption.
\begin{assumption}\label{assumption 4.1} (i)\ The functions $f_i(x)\
(i\in\{0\}\cup I)$ are all second-order
continuously differentiable.\\
(ii)\ The KKT pair $(x^*,\mu^*)$ of problem (\ref{problem1})
satisfies the strong second-order sufficient conditions, i.e.,
$$d^T\nabla^2_{xx}L(x^*,\mu^*)d>0,\ \forall\ d\in R^n,\ d\neq 0,\ g_i(x^*)^Td=0,\ i\in I_*^+,$$
where $\nabla^2_{xx}L(x^*,\mu^*)=\nabla^2f_0(x^*)+\sum\limits_{i\in
I}\mu_i^*\nabla^2f_i(x^*),\ I_*^+=\{i\in I_1:\ \mu_i^*>0\}\cup I_2.$
\end{assumption}
\begin{remark}\label{remark4.2}
Similar to the proof of Lemma \ref{lemma2.3}, we can conclude that
$(x^*, \lambda^*)$ satisfying
 \begin{equation}\label{4.1} \lambda_i^*=\mu_i^*,\
i\in I_1;\ \lambda_i^*=\mu_i^*+c,\ i\in I_2\end{equation}
 is the KKT point of  problem (\ref{problem3}). Moreover,
$\{i\in I:\ \lambda_i^*>0\}=\{i\in I_1:\ \mu_i^*>0\}\cup I_2,$ which
implies that  KKT pair $(x^*, \lambda^*)$ of problem
(\ref{problem3}) also satisfies the strong second-order sufficiency
conditions, i.e.,
$$d^T\nabla^2_{xx}L_c(x^*,\lambda^*)d>0,\ \forall\ d\in R^n,\ d\neq 0,\ g_i(x^*)^Td=0,\ i\in\tilde{I}_*^+,$$
where $\nabla^2_{xx}L_c(x^*,\lambda^*)=\nabla^2F_c(x^*)+
\sum\limits_{i\in I}\lambda_i^*\nabla^2f_i(x^*),\
\tilde{I}_*^+=\{i\in I:\ \lambda_i^*>0\}.$
\end{remark}

Under the stated assumptions, we have the following theorem.
\begin{theorem}\label{theorem4.3}
Suppose that Assumptions \ref{assumption 2.1}, \ref{assumption 3.1}
and \ref{assumption 4.1} hold.
Then\\
(i)$\lim\limits_{k\rightarrow\infty} x^k=x^*,$ i.e., Algorithm
\ref{algorithm 2.7}
is strongly convergent.\\
(ii) $\lim\limits_{k\rightarrow\infty}\lambda^k=\lambda^*$,
$\lim\limits_{k\rightarrow\infty}\mu^k=\mu^*$.
\end{theorem}
{\bf proof:}\ \ \ (i) The proof of this part is similar to the one
of Theorem 4.1 in \cite{jthz2008}, and the details can be seen in
\cite{jthz2008}.

(ii) From the proof of Theorem \ref{theorem3.5} and part (i), one
can conclude that each accumulation point of sequencees
$\{\lambda^k\}$ and $\{\mu^k\}$ is the KKT multiplier for  problem
(\ref{problem3}) and problem (\ref{problem1}) associated with $x^*$
, respectively. Togethering with the uniqueness of the KKT
multiplier, this furthermore implies that part (ii) holds.\qed

\begin{lemma}\label{lemma4.4}
Suppose that Assumptions \ref{assumption 2.1}, \ref{assumption 3.1} and \ref{assumption 4.1} hold. Then\\
(i) $||d_2^k||=O(||d_0^k||^2)+O(\varphi(x^k)^\sigma)$,
$||d_2^k||^2=O(||d_0^k||^4)
+o(\varphi_k^\sigma)$, $||h_2^k||=O(||d_0^k||^2)+O(\varphi(x^k)^\sigma)$.\\
(ii) $\tilde{I}_*^+\subseteq L(x^k)\subseteq I_0(x^*)$ for $k$ large
enough.
\end{lemma}
{\bf proof:}\ \ (i) In view of $\digamma(x^k+d_0^k)=O(||d_0^k||^2)$,
the proof is elementary from (\ref{2.15}) and Lemma
\ref{lemma3.4}(i).

(ii) For $i\not\in I_0(x^*)$, we have $\bar{f}_i(x^*)<0$. Since
$\lim\limits_{k\rightarrow\infty}(x^k, d_0^k)=(x^*, 0)$, there
exists a constant $\bar{\xi}>0$ such that
$\bar{f}_i(x^k)\leq-\bar{\xi}<0$ for $k$ large enough. Moreover, it
 holds that $\bar{f}_i(x^k)+g_i(x^k)^Td_0^k\leq
-\frac{1}{2}\bar{\xi}<0$ for $k$ large enough, which implies that
$i\not\in L(x^k)$, i.e.,  $L(x^k)\subseteq I_0(x^*)$. Furthermore,
it follows from Theorem \ref{theorem4.3}(ii) that
$\lim\limits_{k\rightarrow\infty}\lambda_{\tilde{I}_*^+}^k=\lambda_{\tilde{I}_*^+}^*>0$,
i.e., $\lambda_{\tilde{I}_*^+}^k>0$ and $\tilde{I}_*^+\subseteq
L(x^k)$ hold for $k$ sufficiently large.\qed

It is well-known that the strict complementarity condition (i.e.,
$\tilde{I}^+_*= I_0(x^*)$) is very important to ensure
$L(x^k)=I_0(x^*)$ holds, however, this condition is hard to verify
in practice. In our paper, by the strong second-order sufficient
conditions, we only need $L(x^k)\subseteq I_0(x^*)$.

To ensure the steplength $t_k\equiv1$ for $k$ large enough without
the strict complementary assumption, an additional assumption is
necessary.

\begin{assumption}\label{assumption4.5}  Suppose that the KKT pair
$(x^*,\lambda^*)$ and matrix $H_k$ satisfy
$$||(\nabla^2_{xx}L_c(x^*,\lambda^*)-H_k)d_0^k||=o(||d_0^k||),$$
where
$\nabla^2_{xx}L_c(x^*,\lambda^*)=\nabla^2F_c(x^*)+\sum\limits_{i\in
I}\lambda_i^*\nabla^2f_i(x^*)=\nabla^2L(x^*,\mu^*)$.
\end{assumption}

\begin{theorem}\label{theorem4.6}  Suppose that Assumptions \ref{assumption 2.1},
\ref{assumption 3.1}, \ref{assumption 4.1} and \ref{assumption4.5}
hold. Then the inequalities in (\ref{2.18}) always hold for $t=1$
and $k$ large enough.
\end{theorem}
{\bf proof:}\ \ We assume that $t=1$ and $k$ large enough in the
whole process of proof. First of all, we discuss the second and the
last inequalities of (\ref{2.18}).

For $i\not\in I_0(x^*)$, i.e., $\bar{f}_i(x^*)<0$. In view of
$(x^k,d_0^k,d_2^k,\varphi(x^k))\rightarrow(x^*,0,0,0)$ $
(k\rightarrow\infty)$, we have $d^k=d_0^k+d_2^k\rightarrow0\
(k\rightarrow\infty)$. So we can conclude that the second
inequalities and the last inequality of (\ref{2.18}) are both
satisfied.

For $i\in I_0(x^*)$, it holds that $\bar{f}_i(x^*)=0$. On one hand,
since
$\lim\limits_{k\rightarrow\infty}\bar{f}_i(x^k)=\bar{f}_i(x^*)=0$
and $\lim\limits_{k\rightarrow\infty}\varphi(x^k)=0$ as well as
(\ref{2.12}), it follows that $Q_i^k\rightarrow 0$ and
$Q_i^k=o(|\bar{f}_i(x^k)+g_i(x^k)^Td_0^k|)+o(||d_0^k||).$
 On the
other hand, we have from (\ref{2.15}) and Lemma \ref{lemma4.4}(i)
\begin{equation}\label{4.2}\begin{array}{ll}g_i(x^k)^Td_2^k=-||d_0^k||^\tau-\varphi(x^k)^\sigma-f_i(x^k+d_0^k)+f_i(x^k)+g_i(x^k)^Td_0^k
+o(|\bar{f}_i(x^k)+g_i(x^k)^Td_0^k|)+O(||d_0^k||^3)+o(\varphi(x^k)^\sigma).\end{array}\end{equation}
Then we obtain by Taylor expansion and (\ref{4.2})
\begin{equation}\label{4.3}f_i(x^k+d^k)=\left\{\begin{array}{ll}-||d_0^k||^\tau-|\bar{f}_i(x^k)
+g_i(x^k)^Td_0^k|+o(|\bar{f}_i(x^k)+g_i(x^k)^Td_0^k|)+O(||d_0^k||^3),\
\ \ \ \ \ {\rm{if}}\ i\in
I_0(x^*)\cap I^-;\\
-||d_0^k||^\tau-\varphi(x^k)^\sigma+\varphi(x^k)-|\bar{f}_i(x^k)+g_i(x^k)^Td_0^k|+o(|\bar{f}_i(x^k)+g_i(x^k)^Td_0^k|)\\
\ \ \ +O(||d_0^k||^3)+o(\varphi(x^k)^\sigma)\ \ \ \ \ \ \ {\rm{if}}\
i\in I(x^*)\cap I^+.\end{array}\right.\end{equation} By
$\tau\in(2,3)$, the first equality of (\ref{4.3}) implies that
$f_i(x^k+d^k)\leq 0$ for all $i\in I_0(x^*)\cap I^-$, i.e., the
second inequalities of (\ref{2.18}) hold for $i\in I_0(x^*)\cap I^-$
and the third inequality of (\ref{2.18}) holds.

 Again from (\ref{4.3}), $\tau\in(2,3)$ as well as
$\alpha\in(0,\frac{1}{2})$, we have for $i\in I_0(x^*)\cap I^+$
$$\begin{array}{ll}f_i(x^k+d^k)-\max\{0,\varphi(x^k)-\alpha(||d_0^k||^\tau+\varphi(x^k)^\sigma)\}
&=-(1-\alpha)(||d_0^k||^\tau+\varphi(x^k)^\sigma)-|\bar{f}_i(x^k)+g_i(x^k)^Td_0^k|\\
&\ \ \ \ +o(|\bar{f}_i(x^k)+g_i(x^k)^Td_0^k|)+o(\varphi(x^k)^\sigma)+O(||d_0^k||^3)\\
&\leq 0.\end{array}$$ Summarizing the above analysis, we have proved
that the second and the last inequalities of (\ref{2.18}) are
satisfied for $t=1$ and $k$ large enough.

From now on, we will show that the first inequality of (\ref{2.18})
holds. First of all, by Taylor expansion and Lemma
\ref{lemma4.4}(i), we have
\begin{equation}\label{4.4}\begin{array}{ll}\Delta_k&\triangleq
\nabla F_c(x^k)^Td^k+\frac{1}{2}(d^k)^T\nabla^2F_c(x^k)d^k-\alpha
\nabla F_c(x^k)^Td_0^k-\rho(1-\alpha)\varphi(x^k)^\theta+o(||d^k||^2)\\
&=\nabla
F_c(x^k)^T(d_0^k+d_2^k)+\frac{1}{2}(d_0^k)^T\nabla^2F_c(x^k)d_0^k-\alpha
\nabla F_c(x^k)^Td_0^k-\rho(1-\alpha)\varphi(x^k)^\theta
+o(||d_0^k||^2)+o(\varphi(x^k)^\sigma).\end{array}\end{equation}
Then we get by the KKT conditions (\ref{2.9}) and Lemma
\ref{lemma4.4}(i)
\begin{equation}\label{4.5}\begin{array}{ll} \nabla
F_c(x^k)^T(d_0^k+d_2^k)=-(d_0^k)^TH_kd_0^k-\sum\limits_{i\in
L(x^k)}\lambda_ig_i(x^k)^T(d_0^k+d_2^k)+o(||d_0^k||^2)+o(\varphi(x^k)^\sigma).\end{array}\end{equation}
For $i\in L(x^k)\subseteq I_0(x^*)$, it follows that
$\bar{f}_i(x^k)+g_i(x^k)^Td_0^k=0$. From (\ref{4.3}) and Lemma
\ref{lemma4.4}(i) as well as $\varphi(x^k)=o(\varphi(x^k)^\sigma)$,
we have
$$
f_i(x^k+d^k)=\left\{\begin{array}{ll}
-||d_0^k||^\tau-\varphi(x^k)^\sigma+o(||d_0^k||^2)+o(\varphi(x^k)^\sigma),\
i\in L(x^k);\\
\bar{f}_i(x^k)+g_i(x^k)^T(d_0^k+d_2^k)
+\frac{1}{2}(d_0^k)^T\nabla^2f_i(x^k)d_0^k+o(||d_0^k||^2)+o(\varphi(x^k)^\sigma),\end{array}\right.$$
which further imply that
\begin{equation}\label{4.6}\begin{array}{ll}-\sum\limits_{i\in
L(x^k)}\lambda_i^kg_i(x^k)^T(d_0^k+d_2^k)=\sum\limits_{i\in
L(x^k)}\lambda_i^k\bar{f}_i(x^k)+\frac{1}{2}\sum\limits_{i\in
L(x^k)}\lambda_i^k(d_0^k)^T\nabla^2f_i(x^k)d_0^k+o(||d_0^k||^2)+O(\varphi(x^k)^\sigma).\end{array}\end{equation}
Substituting (\ref{4.6}) into (\ref{4.5}), we have
$$\begin{array}{ll}
\nabla F_c(x^k)^T(d_0^k+d_2^k)=-(d_0^k)^TH_kd_0^k+\sum\limits_{i\in
L(x^k)}\lambda_i^k\bar{f}_i(x^k)+\frac{1}{2}(d_0^k)^T(\nabla^2_{xx}L_c(x^k,\lambda^k)-\nabla^2F_c(x^k))d_0^k+o(||d_0^k||^2)+O(\varphi(x^k)^\sigma),
\end{array}
$$
which combined with (\ref{4.4}) gives
$$\begin{array}{ll}\Delta_k&=-(d_0^k)^TH_kd_0^k+\sum\limits_{i\in L(x^k)}\lambda_i^k\bar{f}_i(x^k)
+\frac{1}{2}(d_0^k)^T\nabla_{xx}^2L_c(x^k,\lambda^k)d_0^k-\alpha
\nabla F_c(x^k)^Td_0^k-\rho(1-\alpha)\varphi(x^k)^\theta+o(||d_0^k||^2)+O(\varphi(x^k)^\sigma)\\
&=(\alpha-\frac{1}{2})(d_0^k)^TH_kd_0^k+(1-\alpha)\sum\limits_{i\in
L(x^k)}\lambda_i^k\bar{f}_i(x^k)+\frac{1}{2}(d_0^k)^T\left(\nabla_{xx}^2L_c(x^k,\lambda^k)-H_k\right)d_0^k\\
&\ \ \ \ -\rho(1-\alpha)\varphi(x^k)^\theta
+o(||d_0^k||^2)+O(\varphi(x^k)^\sigma).\end{array}$$ Due to
$\lambda_i^k\bar{f}_i(x^k)\leq 0$ and $\alpha\in(0,\frac{1}{2})$ as
well as $\theta<\sigma$, it follows from Assumptions \ref{assumption
3.1} and \ref{assumption4.5} that
$$\Delta_k\leq(\alpha-\frac{1}{2})a||d_0^k||^2+o(||d_0^k||^2)-\rho(1-\alpha)
\varphi(x^k)^\theta+o(\varphi(x^k)^\theta)\leq 0,$$ i.e., the first
inequality of (\ref{2.18}) holds. The whole proof is completed.\qed
\begin{theorem}\label{theorem4.7} Under all above mentioned assumptions,
$\varphi(x^{k+1})\equiv0$ after a finite number of iterations, i.e.,
$x^{k+1}\in \Omega^+$ for $k$ large enough.
\end{theorem}

According to Theorem \ref{theorem4.6} and (\ref{4.3}), the above
theorem holds directly. Moreover, based on Theorems \ref{theorem4.6}
and \ref{theorem4.7}, the superlinear convergence of Algorithm
\ref{algorithm 2.7} is given in Theorem \ref{theorem4.8} by Theorem
2.2.3 in \cite{pt1993}.

\begin{theorem}\label{theorem4.8} Suppose that Assumptions
\ref{assumption 2.1}, \ref{assumption 3.1}, \ref{assumption 4.1} and
\ref{assumption4.5} are all satisfied. Then Algorithm \ref{algorithm
2.7} is superlinearly convergent, i.e.,
$||x^{k+1}-x^*||=o(||x^k-x^*||)$.
\end{theorem}

\section{Numerical experiments}\label{sec5}
In this section, in order to show the effectiveness of our proposed
algorithm, some classical problems in \cite{hs1980,got2003} are
tested and the corresponding comparative numerical results are
reported in the following parts. The algorithm is implemented by
using MATLAB R2008a on Windows XP platform, and on a PC with 2.53GHz
CPU.

During the numerical experiments, the identity matrix $E_n$ is
selected as the initial Lagrangian Hessian, and the approximation
Hessian matrix $H_k$ is updated by BFGS formula \cite{pm1991}
$$H_{k+1}=H_k-\frac{H_ks^k(s^k)^TH_k}{(s^k)^TH_ks^k}+\frac{\hat{y}^k(\hat{y}^k)^T}{(s^k)^T\hat{y}^k}\
\ (k\geq 0),$$
 where
$$\left\{\begin{array}{llll}s^k=x^{k+1}-x^k,\ \hat{y}^k=y^k+\alpha_k(\gamma_ks^k+A_kA_k^Ts^k),\
\gamma_k=\min\{||d_0^k||^2,\ \kappa\in(0,1)\},\ A_k=( g_i(x^k), i\in L(x^k)),\\
y^k=\nabla_x L_{c_k}(x^{k+1},\lambda_k)-\nabla_x L_{c_k}(x^k,
\lambda_k),\ \nabla_x L_{c_k}(x^k, \lambda^k)=\nabla
F_{c_k}(x^k)+\sum\limits_{i\in
I}\lambda_i^kg_i(x^k),\end{array}\right.$$ and
$$\alpha_k=\left\{\begin{array}{ll}0, &{\rm{if}}\ (s^k)^Ty^k\geq
\mu||s^k||^2,\ \mu\in(0,1),\\ 1, &{\rm{if}}\ 0\leq
(s^k)^Ty^k<\mu||s^k||^2,\\
1+\frac{\gamma_k||s^k||^2-(s^k)^Ty^k}{\gamma_k||s^k||^2+(s^k)^TA_k(A_k)^Ts^k},\
&{\rm{otherwise}}.\end{array}\right.$$ The parameters are selected
as follows
$$\begin{array}{ll}\alpha=\hat{\alpha}=\eta=\kappa=\mu=c_{-1}=0.5,\ \theta=0.4,\
\sigma=0.6,\ \rho=2,\ \tau=2.5,\ \epsilon=0.5^3,\ p=2,\ \gamma_0=2,\
\gamma=1.\end{array}$$ The algorithm stops if the termination
criterions $\parallel d_0^k\parallel<=\bar{\varepsilon}$ and
$\varphi(x)=0$ are both satisfied.

First of all, some notations used in the following tables are
defined in Table \ref{tab: Notation in tables of results}.
\begin{table}[htbp] \centering \footnotesize\caption{Definitions of some notations}
\begin{tabular}{ ll|lllllllllllllllllllllllll }
 \hline\hline
Prob         & The number of test problem in \cite{hs1980}. &Nf0 & The number of objective function evaluations.\\
$n$          & The number of variables of test problem.     &Nf & The number of all constraint functions evaluations.\\
$I_1/I_2$ &The number of equality and inequality
constraints,respectively.                                   &Fv &The objective function value at the final iteration point.\\
Nio          & The number of iterations out of the feasible set.&CPU          & The CPU time (second).\\
Nii & The number of iterations within the feasible set. &$-$&The information is not given in the corresponding references.\\
Ni &The total number of iterations, i.e., Ni=Nio+Nii.&&\\
\hline\hline
\end{tabular}
\label{tab: Notation in tables of results}
\end{table}

In order to show the computational efficiency of Algorithm
\ref{algorithm 2.7} (shorted by ALGO \ref{algorithm 2.7}), which is
compared with other types of algorithms, including SQP algorithms
and systems of linear equations (SLE) algorithms. The statistics of
these algorithms are given in Table \ref{tab: Statistics for
compared algorithm}. The ``Feasible" (or ``Infeasible") in Table
\ref{tab: Statistics for compared algorithm} means that the initial
iteration point have to be feasible (or can be chosen arbitrarily)
for the solving problem.
\begin{table}[htbp] \centering \footnotesize\caption{Description of some comparative algorithms}
\begin{tabular}{ lllllllllllllllllllllllllll}
 \hline\hline
  Author                  &Types of algorithm (shorted by) &Types of solving problem  \\
 \hline
 Jin and Wang \cite{jw2010}  &Feasible SQP (JW-FSQP)                    &nonlinear inequality constrained programming\\
 Qi and Yang \cite{qy2001}   &Infeasible SQP (QY-IFSQP)                    &nonlinear inequality constrained programming\\
 Wang, Chen and He \cite{wch2005}&Infeasible SLE (WCH-IFSLE) &nonlinear equality and inequality constrained programming\\
 Guo \cite{g2011} &Infeasible SQP (G-IFSQP) &nonlinear inequality constrained programming\\
 Jian, Ke, Zheng and Tang \cite{jkzt2009}&Infeasible SQP (JKZT-IFSQP)&nonlinear inequality constrained programming\\
 Yang, Li and Qi \cite{ylq2003}&Feasible SLE (YLQ-FSLE)&nonlinear inequality constrained
 programming\\
 Gu and Zhu \cite{gz2011}&Infeasible SQP (GZ-IFSQP)&nonlinear equality and inequality constrained programming\\
 Gill, Murray and Saunders \cite{gms2005}&SNOPT &nonlinear equality and inequality  constrained programming\\
 \hline\hline
\end{tabular}
\label{tab: Statistics for compared algorithm}
\end{table}

 In Table \ref{tab: Numerical results for
ALGO 2.1 and JW-FSQP},  we compare the number of Ni and Fv required
by ALGO
 \ref{algorithm 2.7} with those required by JW-FSQP. The test problems
 are chosen from \cite{hs1980}, and initial iteration points are all feasible except Prob 030.
 The optimality
tolerance is the same as in \cite{jw2010}. The results in Table
\ref{tab: Numerical results for ALGO 2.1 and JW-FSQP} show that the
number of iterations of ALGO \ref{algorithm 2.7} is much smaller
than that of JW-FSQP for most test problems. From the viewpoints of
Ni and Fv, we can conclude that ALGO \ref{algorithm 2.7} is more
effective than JW-FSQP.
\begin{table}[htbp] \centering \footnotesize\caption{Comparative numerical results of ALGO 2.1 and JW-FSQP}
\begin{tabular}{ lllllllllllllllllllllllllll }
 \hline\hline
   Prob     &$n/|I_1|/|I_2|$      &ALGO 2.1&&&     &JW-FSQP &&\\
          \cline{3-5}\cline{7-9}  & &Ni          &Fv   &CPU     &     &Ni              &Fv\\
 \hline
   012      &2/1/0            &7           &$-3.000000000E+01$   &0.02 &   &18  &$-3.000000000E+01$ \\
 \hline
  024       &2/5/0            &11            &$-1.000000000E+00$   &0.06 &   &12  &$-9.999997439E-01$ \\
  \hline
   029      &3/1/0            &10           &$-2.262741700E+01$   &0.05 &   &16  &$-2.262741700E+01$   \\
 \hline
  030       &3/7/0            &10            &$4.016837909E-18$   &0.05 &   &14  &$1.000000000E+00$ \\
 \hline
   031     &3/7/0             &13          &$6.000000000E+00$   &0.08 &   &12  &$6.000000007E+00$   \\
 \hline
  033       &3/6/0            &9          &$-4.585785958E+00$   &0.05 &   &44  &$-4.585786409E+00$ \\
  \hline
   043     &4/3/0            &10           &$-4.400000000E+01$   &0.05 &   &24  &$-4.400000000E+01$   \\
 \hline
  076      &4/7/0            &9           &$-4.681818182E+00$   &0.05 &  &9  &$-4.681818182E+00$ \\
 \hline
   100      &7/4/0           &15            &$6.825663841E+02$   &0.11&    &42  &$6.806300574E+02$ \\
\hline\hline
\end{tabular}
\label{tab: Numerical results for ALGO 2.1 and JW-FSQP}
\end{table}
\begin{table}[htbp] \centering \footnotesize\caption{Comparative numerical results of ALGO 2.1 and QY-IFSQP}
\begin{tabular}{ lllllllllllllllllllllllllll }
 \hline\hline
   Prob   &$n/|I_1|/|I_2|$  &Point &                          &Ni  &Nf0  &Nf  &Fv                   &CPU\\
 \hline
   034    &3/8/0           &(a)    &ALGO 2.1     &24  &25   &443 &$-0.83403244521568$ &$0.17$\\
          &  &                     &QY-IFSQP                        &34  &161  &162 &$-0.83403244524796$ &$-$\\
          &               &(b)     &ALGO 2.1     &5+12 &18 &374 &$-0.83403244522367$&$0.13$\\
          &               &        &QY-IFSQP                        &13  &32   &33 &$-0.83403244526530$&$-$\\
 \hline
   035    &3/4/0           &(a)    &ALGO 2.1     &12  &13   &119 &$0.11111111111111$ &$0.05$\\
          &  &                     &QY-IFSQP                        &10  &46  &46 &$0.11111111111111$ &$-$\\
          &               &(b)     &ALGO 2.1     &1+9 &11 &84 &$0.11111111111111$&$0.03$\\
          &               &        &QY-IFSQP                       &9&11&12&$0.11111111111111$&$-$\\
 \hline
   036    &3/7/0           &(a)    &ALGO 2.1     &7  &8   &114 &$-3299.99999999996$ &$0.02$\\
          &  &                     &QY-IFSQP                        &3  &7  &7 &$-3300$ &$-$\\
          &               &(b)     &ALGO 2.1      &1+4 &6 &77 &$-3299.99999999997$&$0.03$\\
          &               &        &QY-IFSQP                       &5&10&10&$-3300$&$-$\\
 \hline
   037    &3/8/0           &(a)    &ALGO 2.1     &23  &24   &467 &$-3455.999999999965$ &$0.17$\\
          &  &                     &QY-IFSQP                        &12  &34  &34 &$-3456.000000000001$ &$-$\\
          &               &(b)     &ALGO 2.1     &34+35&70&1189&$-3455.999999999998$&$0.58$\\
          &               &        &QY-IFSQP                       &29&113&119&$-3456.000000000001$&$-$\\
 \hline
   043    &4/3/0           &(a)    &ALGO 2.1     &12  &13   &81 &$-44$ &$0.06$\\
          &  &                     &QY-IFSQP                        &14  &29  &29 &$-44.00000000000001$ &$-$\\
          &               &(b)     &ALGO 2.1     &67+8&76&1270&$-44$&$0.56$\\
          &               &        &QY-IFSQP                       &19&60&65&$-44.00000000000011$&$-$\\
 \hline
   044    &4/10/0           &(a)    &ALGO 2.1     &20  &21   &448 &$-14.99999999935652$ &$0.17$\\
          &  &                     &QY-IFSQP                        &7  &13  &13 &$-15$ &$-$\\
          &               &(b)     &ALGO 2.1     &4+5&10&190&$-14.99999999999756$&$0.09$\\
          &               &        &QY-IFSQP                               &8&18&21&$-15$&$-$\\
 \hline
   065    &3/7/0           &(a)    &ALGO 2.1     &8  &9   &126 &$0.95352885680478$ &$0.03$\\
          &  &                     &QY-IFSQP                        &13  &38  &39 &$0.95352885680478$ &$-$\\
          &               &(b)     &ALGO 2.1      &1+13&15&217&$0.95352885680478$&$0.11$\\
          &               &        &QY-IFSQP                        &15&41&44&$0.95352885680188$&$-$\\
 \hline
   066    &3/8/0           &(a)    &ALGO 2.1     &10  &11   &175 &$0.51816327418156$ &$0.06$\\
          &  &                     &QY-IFSQP                        &8  &8  &8 &$0.51816327418154$ &$-$\\
          &               &(b)     &ALGO 2.1     &2+13&16&252&$0.51816327418154$&$0.11$\\
          &               &        &QY-IFSQP                               &14&26&28&$0.51816327418153$&$-$\\
 \hline
   100    &7/4/0           &(a)    &ALGO 2.1     &20  &21   &193 &$682.5663838261504$ &$0.14$\\
          &  &                     &QY-IFSQP                       &28  &96  &96 &$680.6300573744018$ &$-$\\
          &               &(b)     &ALGO 2.1      &7+14&22&258&$682.5663838261520$&$0.16$\\
          &               &        &QY-IFSQP                               &26&90&92&$680.6300573743961$&$-$\\
\hline\hline
\end{tabular}
\label{tab: Numerical results for ALGO 2.1 and QY-SQP}
\end{table}

In Table \ref{tab: Numerical results for ALGO 2.1 and QY-SQP}, we
further compare the number of Ni, Nf0, Nf and Fv required by ALGO
\ref{algorithm 2.7} with those required by QY-IFSQP. The optimality
criterions and the starting iteration points for the test problems
are the same as in \cite{qy2001}. ``Point (a)" (or ``Point (b)") in
Table \ref{tab: Numerical results for ALGO 2.1 and QY-SQP} denotes
that the corresponding initial point is ``feasible" (or
``infeasible"). Note that Nf refers to the number of evaluations of
$f_I$ in \cite{qy2001}, however, Nf denotes the number of all
constraint functions evaluations in our paper. The Ni column is
displayed as the total number of iterations. Only if the initial
iteration points are chosen as (a), Ni=Nii, otherwise, Ni=Nio+Nii.
For example, ``5+12" means that the algorithm generates a feasible
point after five iterations, and after another twelve iterations the
algorithm produces an approximately optimal solution. For the test
problems, from the viewpoints of Ni and Nf0, the results show that
ALGO \ref{algorithm 2.7} is obviously better than QY-IFSQP for Prob
034, 043, 065 and 100 for point (a). For point (b), again from the
viewpoints of Ni, Nf0 and Fv, ALG \ref{algorithm 2.7} is competitive
with QY-IFSQP for most of test problems except Prob. 037 and 043.

Note that all test problems in Tables \ref{tab: Numerical results
for ALGO 2.1 and JW-FSQP} and \ref{tab: Numerical results for ALGO
2.1 and QY-SQP} only have inequality constraints. In order to show
the performance of ALGO \ref{algorithm 2.7} for solving problems
with equality constraints, ALGO \ref{algorithm 2.7} is further
compared with WCH-IFSLE, GZ-IFSQP and SNOPT, respectively. The test
problems and stopping criterions as well as initial iteration points
are the same as in \cite{wch2005} and \cite{gz2011}, respectively.

For comparing the performance of ALGO \ref{algorithm 2.7} with
WCH-IFSLE and GZ-IFSQP as well as SNOPT, we use performance profiles
as described in Dolan and Mor\'{e}'s paper \cite{dm2002}. Our
profiles for figures are based on the number of iterations. The
function $\rho(\tau)$ is the (cumulative) distribution function for
the performance ratio within a factor $\tau\in R$. The value of
$\rho(\tau)$ is the probability that the solver will win over the
rest of the solvers. The corresponding results of performance are
shown in Figure \ref{fig:Performance profiles for ALGO 2.1 and
WCH-SLE}. From Figure \ref{fig:Performance profiles for ALGO 2.1 and
WCH-SLE}, it is obviously that the performance of ALGO
\ref{algorithm 2.7} is better than that of WCH-IFSLE, i.e., ALGO
\ref{algorithm 2.7} has the most wins compare with WCH-SLE.
Moreover, our algorithm is competitive with SNOPT (which is a
well-known SQP algorithm for solving nonlinear constrained
programming) although the performance of GZ-IFSQP is better than
ALGO \ref{algorithm 2.7}.
\begin{figure}[htbp]\centering
\includegraphics[width=8cm]{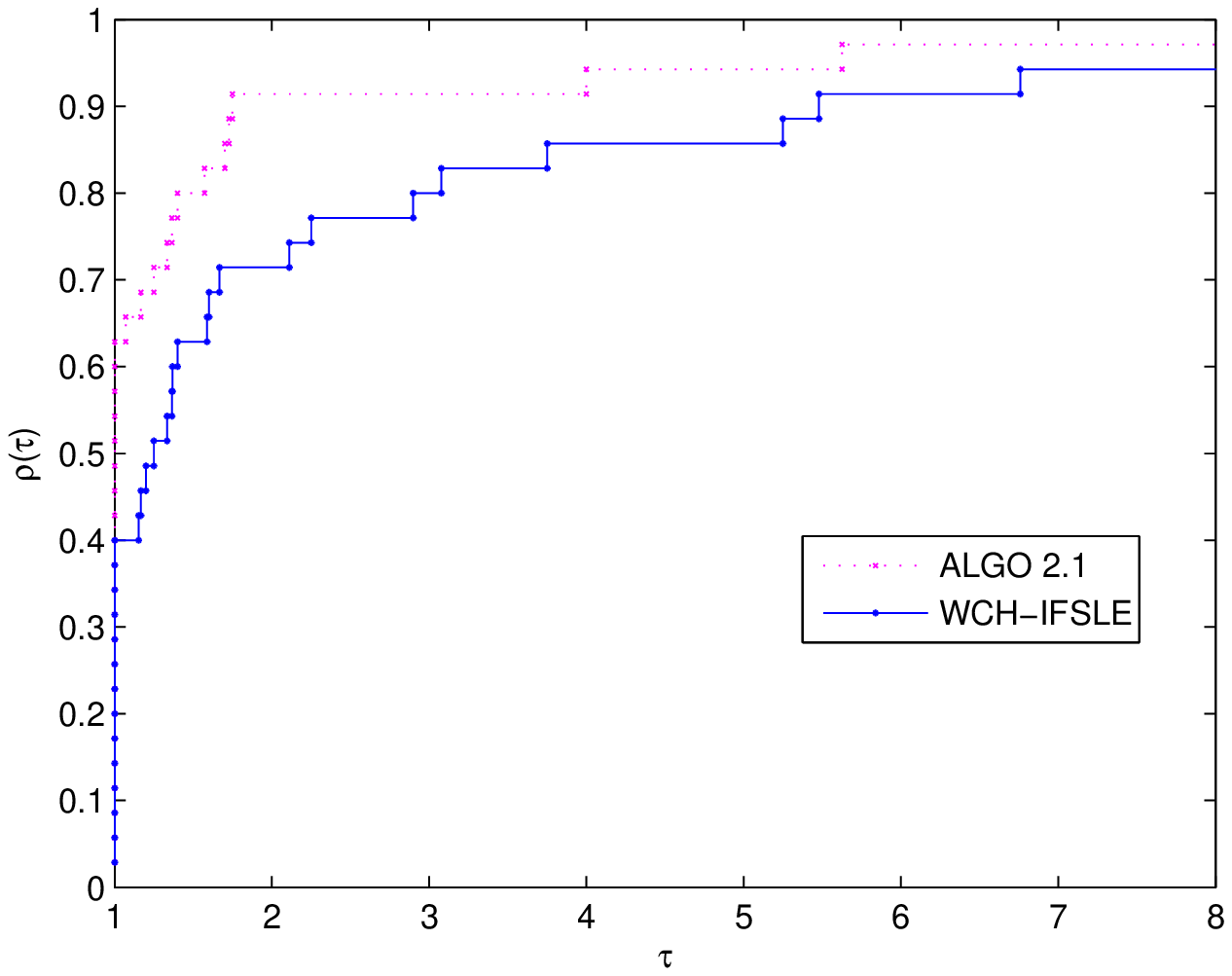}
\includegraphics[width=8cm]{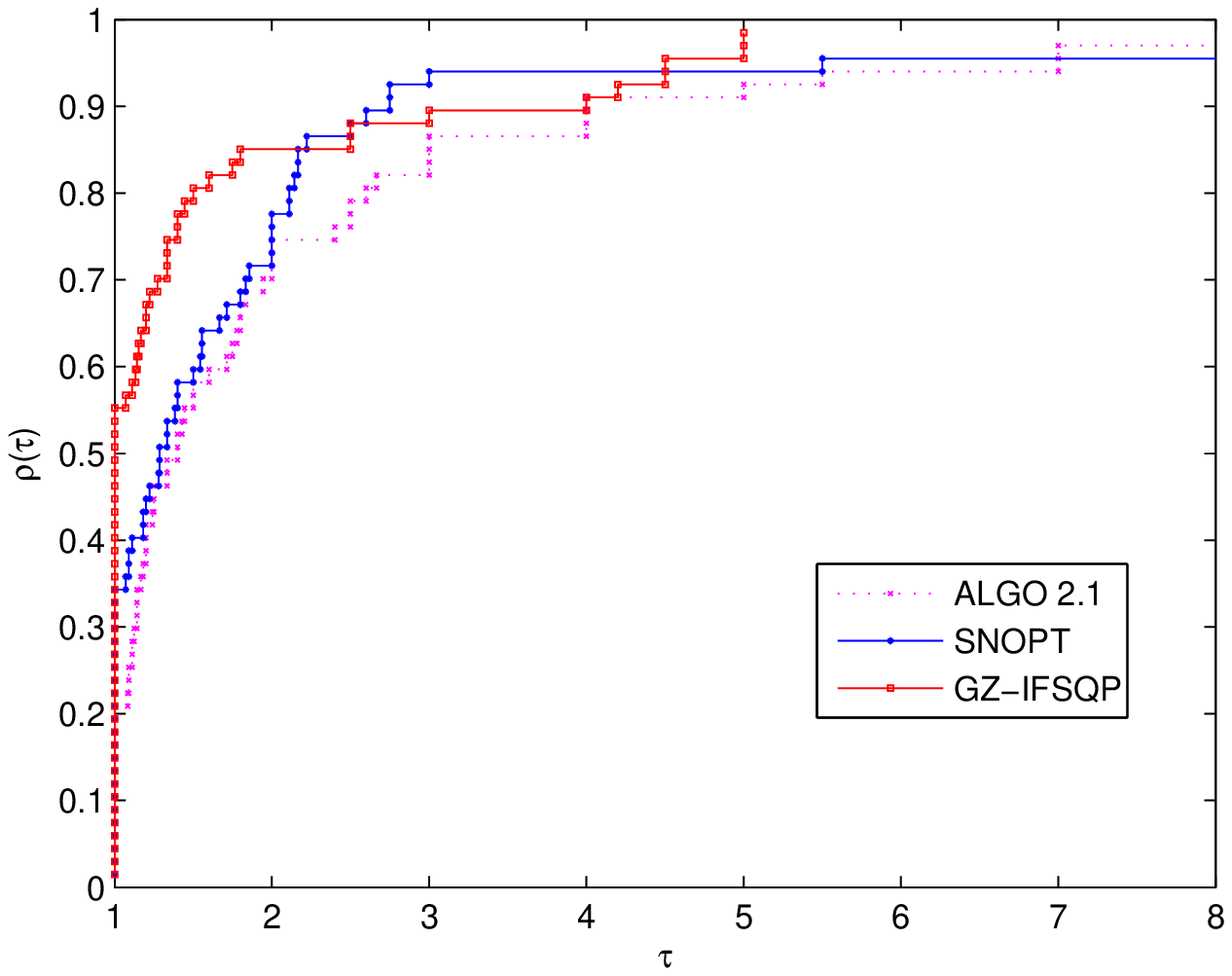}
\footnotesize\caption{The left figure shows the performance of ALGO
2.1 and WCH-IFSLE, the right figure shows the performance of ALGO
2.1 and GZ-IFSQP as well as SNOPT.}\label{fig:Performance profiles
for ALGO 2.1 and WCH-SLE}
\end{figure}
\begin{table}[htbp] \centering \footnotesize\caption{Comparative numerical results of ALGO 2.1, G-IFSQP, JKZT-IFSQP and YLQ-FSLE}
\begin{tabular}{ lllllllllllllllllllllllllll }
 \hline\hline
   Prob         &$n/|I_1|/|I_2|$   &Algorithm                    &Nii  &Nf0  &Nf      &Fv          &CPU\\
\hline
   Svanberg-10  &10/30/0           &ALGO 2.1 &14   &15    &1140      &15.731517   &0.38\\
                   &&G-IFSQP &17 &18 &1200 &15.731517 &0.44\\
                  &&JKZT-IFSQP &28 &28 &1753 &15.731533 &$-$\\
                   &&YLQ-FSLE &36 &227 &258 &15.731517 &$-$\\
\hline
   Svanberg-30  &30/90/0           &ALGO 2.1 &23   &24    &6210      &49.142526  &2.67\\
                   &&G-IFSQP &25 &26 &5490 &49.142526 &2.09\\
                  &&JKZT-IFSQP &27 &27 &4975 &49.142545 &$-$\\
                   &&YLQ-FSLE &101 &777 &864 &49.142526 &$-$\\
\hline
   Svanberg-50  &50/150/0           &ALGO 2.1 &29   &30    &13050      &82.581912   &5.59\\
                   &&G-IFSQP &33 &34 &11550 &82.581912 &5.95\\
                  &&JKZT-IFSQP &37 &37 &11762 &82.581928 &$-$\\
                   &&YLQ-FSLE &108 &881 &968 &82.581912 &$-$\\
\hline
   Svanberg-80  &80/240/0           &ALGO 2.1 &38   &39    &33120      &132.749819   &13.80\\
                   &&G-IFSQP &42 &43 &24720 &132.749819 &15.38\\
                  &&JKZT-IFSQP &47 &47 &24100 &132.749830 &$-$\\
                   &&YLQ-FSLE &190 &1666 &1835 &132.749819 &$-$\\
\hline
   Svanberg-100  &100/300/0           &ALGO 2.1 &42   &43    &43200     &166.197171   &23.09\\
                   &&G-IFSQP &55 &56 &39600 &166.197171 &30.53\\
                  &&JKZT-IFSQP &46 &46 &27880 &166.197199 &$-$\\
                   &&YLQ-FSLE &178 &1628 &1782 &166.197171 &$-$\\
 \hline\hline
\end{tabular}
\label{tab: Numerical results for ALGO 2.1, G-SQP, JKZT-SQP and
YLQ-SLE}
\end{table}

Note that the above test problems are relatively small. In order to
show the more clearly effectiveness of ALGO \ref{algorithm 2.7} for
solving some large scale problems, the ``Svanberg" problems are
tested, which are selected from CUTE \cite{got2003}. The
corresponding results are given in Tables \ref{tab: Numerical
results for ALGO 2.1, G-SQP, JKZT-SQP and YLQ-SLE} and \ref{tab:
Numerical results for ALGO 2.1 and G-SQP}. In Table \ref{tab:
Numerical results for ALGO 2.1, G-SQP, JKZT-SQP and YLQ-SLE}, the
performance of ALGO \ref{algorithm 2.7} is compared with G-IFSQP,
JKZT-IFSQP and YLQ-FSLE, respectively. The initial iteration points
are feasible and the stopping criterions are the same as that
reported in \cite{g2011}. From the results in Table \ref{tab:
Numerical results for ALGO 2.1, G-SQP, JKZT-SQP and YLQ-SLE}, in
viewpoints of NII and NF0, it follows that ALGO \ref{algorithm 2.7}
is more effective than G-IFSQP and JKZT-SQP as well as YLQ-SLE for
solving ``Svanberg" problems, respectively.

\begin{table}[htbp] \centering \footnotesize\caption{Comparative numerical results of ALGO 2.1 and G-IFSQP}
\begin{tabular}{ lllllllllllllllllllllllllll }
 \hline\hline
   Prob         &$n/|I_1|/|I_2|$   &Algorithm                  &Nio  &Nii  &Nf0  &Nf      &Fv          &CPU\\
\hline
   Svanberg-10  &10/30/0           &ALGO 2.1  &3   &15  &19  &1290      &15.731517   &0.42\\
                  &&G-IFSQP &3 &15&19 &1320 &15.731517 &0.73\\
\hline
   Svanberg-20  &20/60/0           &ALGO 2.1 &4   &21&26    &4560      &32.427932  &1.53\\
                  &&G-IFSQP &4&22 &27 &4380 &32.427932 &1.69\\
\hline
   Svanberg-30  &30/90/0           &ALGO 2.1 &3   &23&27    &6570      &49.142526   &3.06\\
                  &&G-IFSQP &3 &25&29 &6480 &49.142526 &2.58\\
\hline
   Svanberg-40  &40/120/0           &ALGO 2.1 &3   &23&27    &9840      &65.861140   &3.75\\
                  &&G-IFSQP &3&28 &32 &9480 &65.861140 &3.83\\
\hline
   Svanberg-50  &50/150/0           &ALGO 2.1 &5   &36&42    &23700     &82.581912   &8.27\\
                  &&G-IFSQP &11&30 &42 &16050 &82.581912 &8.16\\
\hline
   Svanberg-80  &80/240/0           &ALGO 2.1 &5   &81&87    &101280     &132.749819   &34.11\\
                  &&G-IFSQP &2&48 &51 &28320 &132.749819 &18.42\\
\hline
   Svanberg-100  &100/300/0           &ALGO 2.1 &3   &43 &47   &43500     &166.197171   &25.58\\
                  &&G-IFSQP &3&52 &56 &39900 &166.197171 &31.05\\
\hline
   Svanberg-150  &150/450/0           &ALGO 2.1 &40   &44&85    &166500     &249.818369   &112.69\\
                  &&G-IFSQP &34&59 &94 &153000 &249.818369 &127.84\\
\hline
   Svanberg-200  &200/600/0           &ALGO 2.1 &3   &84&88    &232800     &333.441310   &226.14\\
                  &&G-IFSQP &4&94 &99 &148200 &333.441310 &307.72\\
 \hline\hline
\end{tabular}
\label{tab: Numerical results for ALGO 2.1 and G-SQP}
\end{table}
Moreover, ALGO \ref{algorithm 2.7} is compared with G-IFSQP
  for ``Svanberg" problems with infeasible initial iteration point, i.e. $x^0=(10,\ldots,10)^T$. The
  optimality thresholds are the same as in
  \cite{g2011}, and the comparative results are given in Table \ref{tab: Numerical results for ALGO 2.1 and G-SQP}. For
  Prob Svanberg-20, Svanberg-30 and Svanberg-40, although Nio and Fv of ALGO \ref{algorithm 2.7} are the
  same as that of G-IFSQP, Nii of ALGO \ref{algorithm 2.7}
  is less than that of G-IFSQP. For Prob Svanberg-50 and Svanberg-200,
  ALGO \ref{algorithm 2.7} can enter into the feasible region more
  quickly than G-SQP. In view of Nio, Nii and Fv, it holds that
  ALGO \ref{algorithm 2.7} is more competitive than G-IFSQP.
  Furthermore, the comparative results of ALGO \ref{algorithm 2.7} and
G-IFSQP in Tables \ref{tab: Numerical results for ALGO 2.1, G-SQP,
JKZT-SQP and YLQ-SLE} and \ref{tab: Numerical results for ALGO 2.1
and G-SQP} further imply that the efficiency of MQSSFD is higher
than that of MSSFD. And this also proves that the conclusions in
\cite{jcg2011} is correct.

\section{Conclusions}\label{sec6}
In this paper, inspired by the ideas in \cite{mp1976,g2011}, an
improved SQP algorithm with arbitrary initial iteration point for
  solving problem (\ref{problem1}) is proposed. Firstly, problem (\ref{problem1}) is equivalently transformed
  into an associated simpler problem (\ref{problem2}). At each iteration, the search direction is generated by solving
  an always QP subproblem and one (or two) SLE (s). The two SLEs have the same coefficient
  matrices. After a finite number of iterations, the iteration points always lie in the feasible region of
   problem (\ref{problem2}), and we only need to
solve the one SLE.
  In the process of iteration, the feasibility of the iteration
 points is monotone increasing. Under some mild assumptions without the strict complementary,
 our algorithm possesses global and superlinear convergence.
 Some comparative numerical results in Section \ref{sec5} show that our algorithm is effective and promising.

\section*{Acknowledgements}The authors would like to thank the associated
editor and the one anonymous referee for taking the time to provide
detailed and highly valuable comments, which significantly improved
the quality of our manuscript. The first author would also like to
thank Dr. Li-Ping Tang and Dr. Jing Zhang for their help to revise
English language errors in the manuscript.


\bibliographystyle{model1-num-names}
\bibliography{references_abbrev}


\end{document}